\documentclass{aims}
\usepackage{amssymb,amsmath,mathtools}
  \usepackage{paralist}
  \usepackage{graphics} 
  \usepackage{epsfig} 
\usepackage{graphicx}  \usepackage{epstopdf}
 \usepackage[colorlinks=true]{hyperref}
\hypersetup{urlcolor=blue, citecolor=red}

\def\currentvolume{X}
 \def\currentissue{X}
  \def\currentyear{200X}
   \def\currentmonth{XX}
    \def\ppages{X--XX}

\theoremstyle{definition}
\newtheorem*{aww}{A word of warning}
\theoremstyle{remark}
\newtheorem*{remark}{Remark}
\def\R {\mathbb{R}}
\def\A {\mathbb{A}}

\def\E {{\mathsf{E}}}

\def\FigA{\begin{figure}[htb]\begin{center}
\includegraphics[width=7cm]{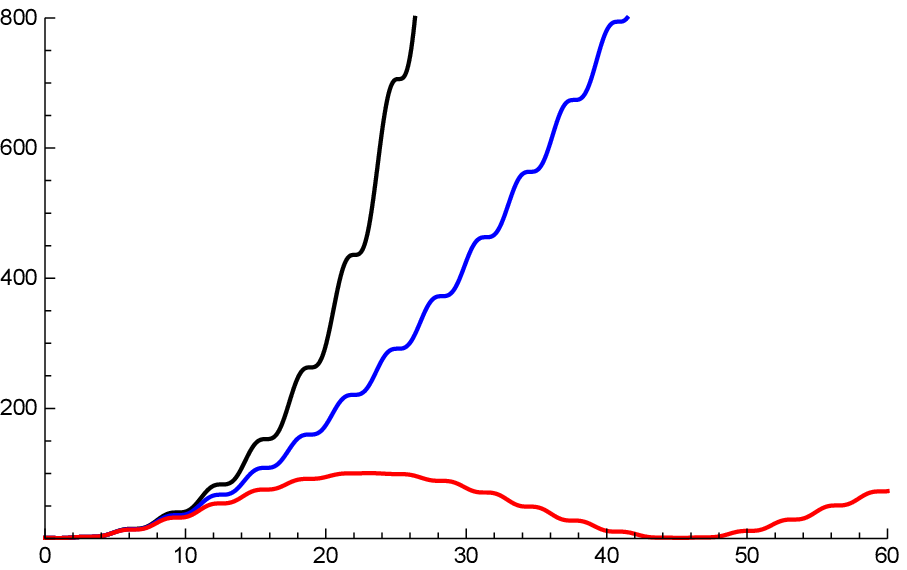}\\
{\tiny Fig.\ $\!$1$\,\,$ Plot of $\E$ for $\epsilon=1$ and $b=0.99$ (black), $b=1$ (blue) and $b=1.01$ (red).}
\end{center}
\end{figure}}

\def\FigB{\begin{figure}[htb]\begin{center}
\includegraphics[width=7cm]{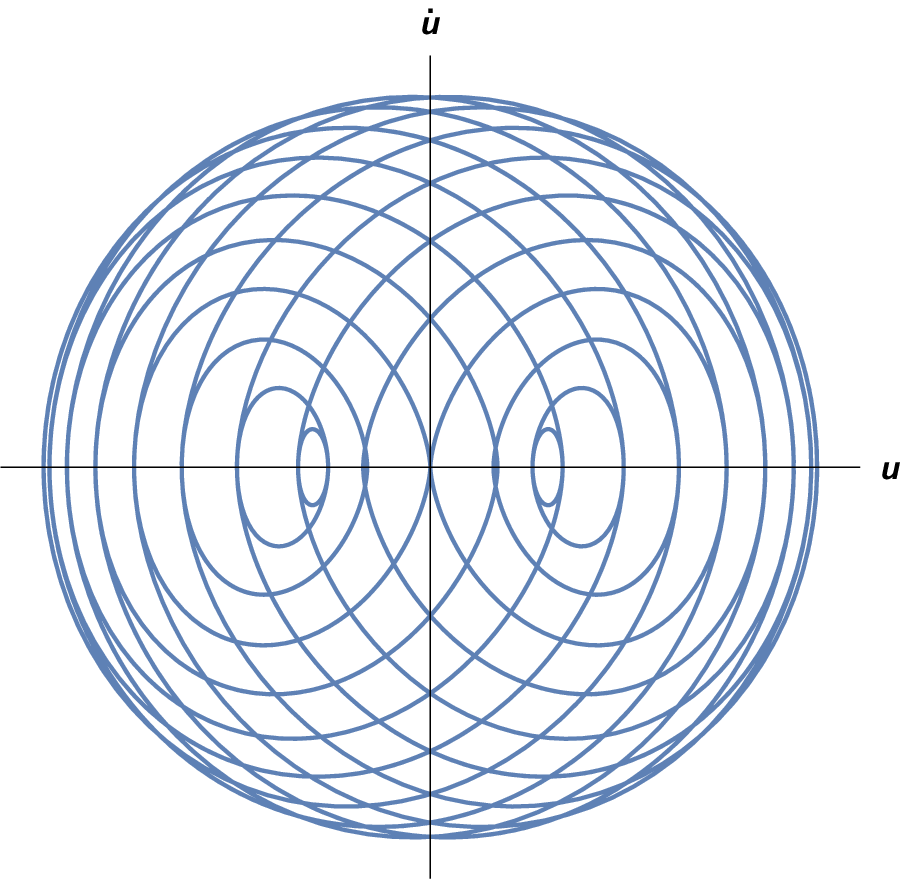}\\
{\tiny Fig.\ $\!$2$\,\,$ Parametric plot of $t\mapsto (u(t),\dot  u(t))$
for $\epsilon=1$ and $b=\sqrt{\frac{23}{13}+\frac{13}{23}-1}$.}
\end{center}
\end{figure}}

\def\FigC{\begin{figure}[htb]\begin{center}
\includegraphics[width=4cm]{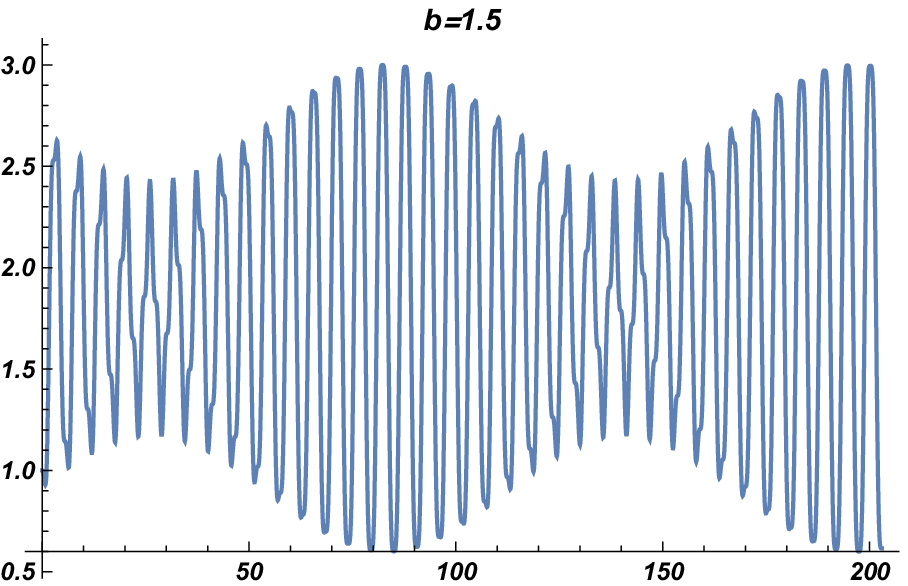}\,\,\,\includegraphics[width=4cm]{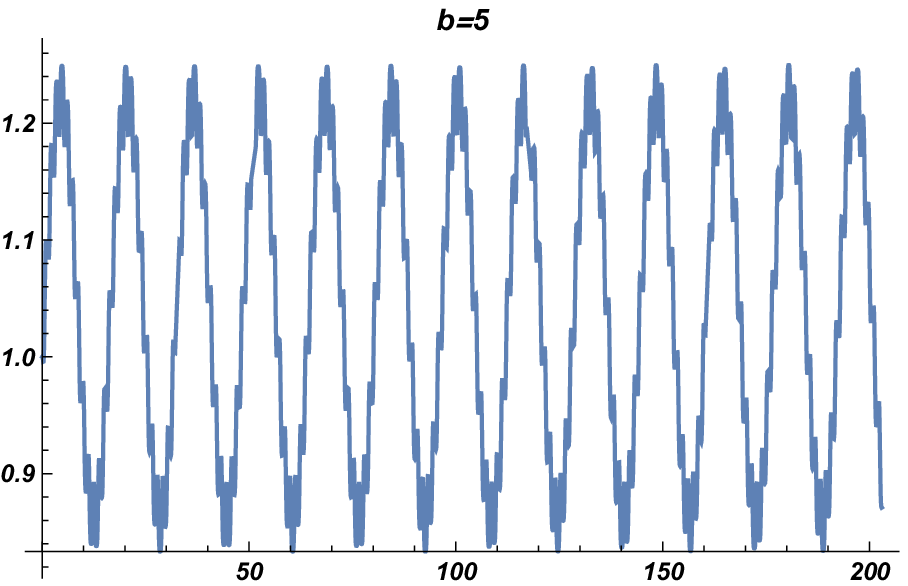}\,\,\,{\includegraphics[width=4cm]{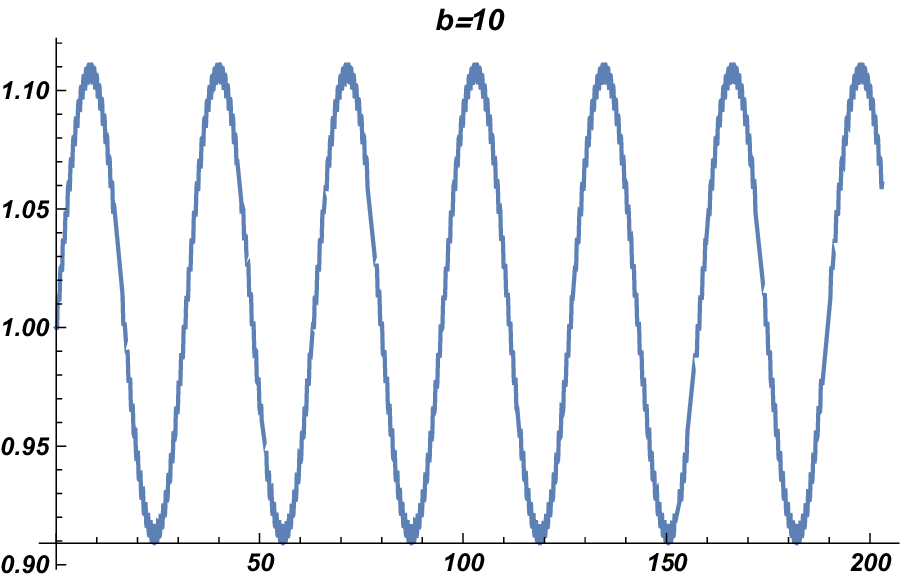}}\\
{\tiny Fig.\ $\!$3$\,\,$ $\E$ for $\epsilon=1$ and $\boldsymbol{z}_0=(1,0,0,0)$
with different values of $b$.}
\end{center}
\end{figure}}

\def\FigD{\begin{figure}[htb]\begin{center}
\includegraphics[width=4cm]{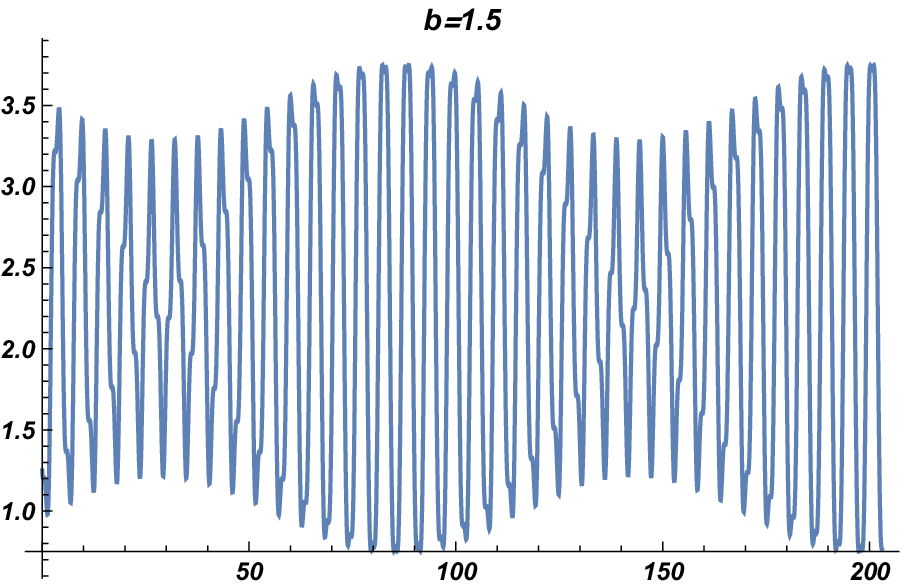}\,\,\,\includegraphics[width=4cm]{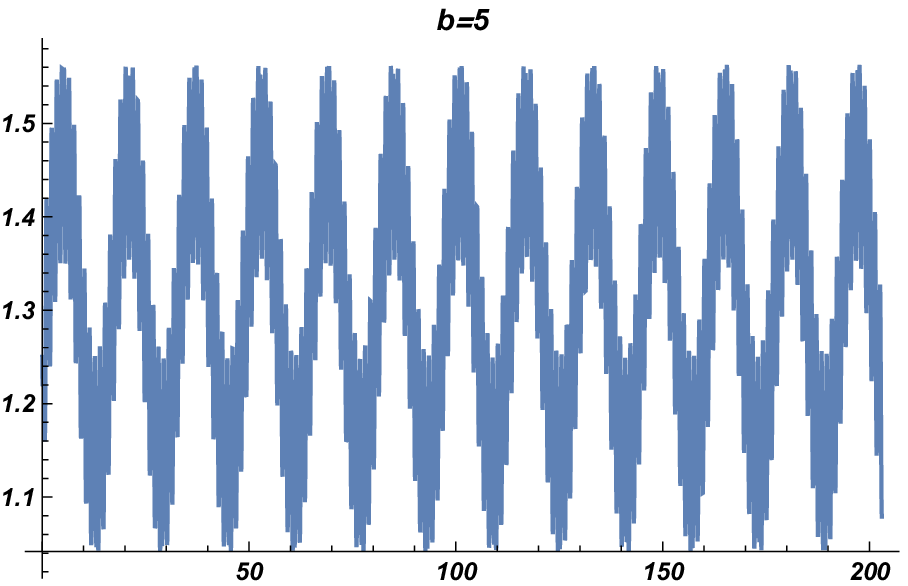}\,\,\,{\includegraphics[width=4cm]{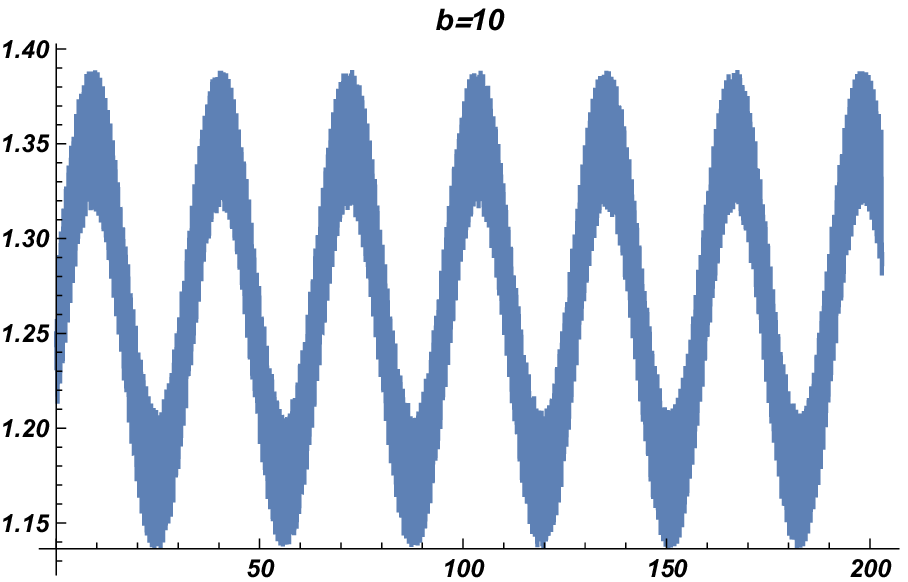}}\\
{\tiny Fig.\ $\!$4$\,\,$ $\E$ for $\epsilon=1$ and $\boldsymbol{z}_0=(1,0.5,0,0)$
with different values of $b$.}
\end{center}
\end{figure}}

\def\FigE{\begin{figure}[htb]\begin{center}
\includegraphics[width=4cm]{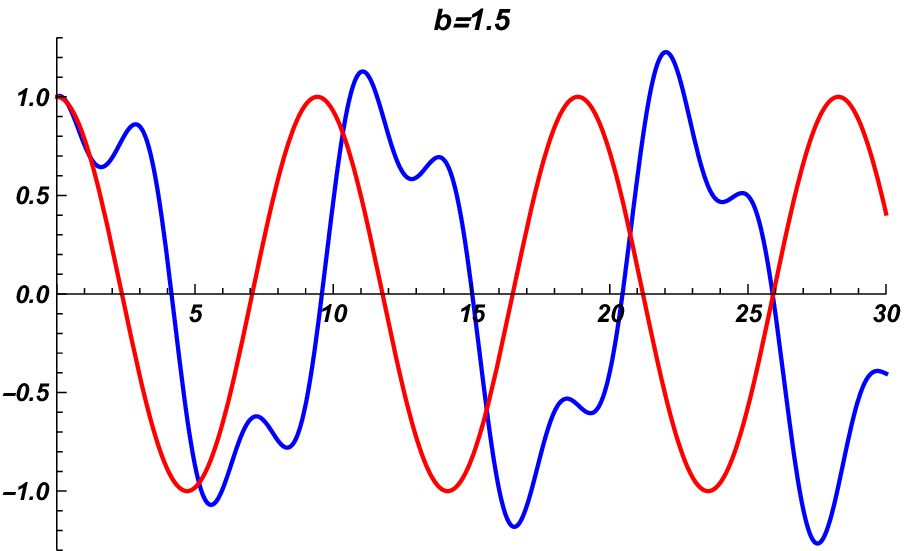}\,\,\,\includegraphics[width=4cm]{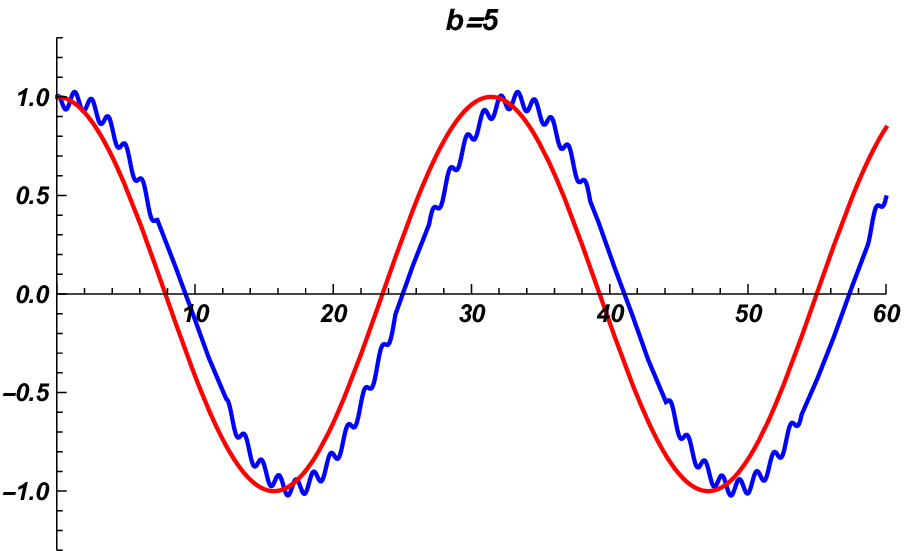}\,\,\,{\includegraphics[width=4cm]{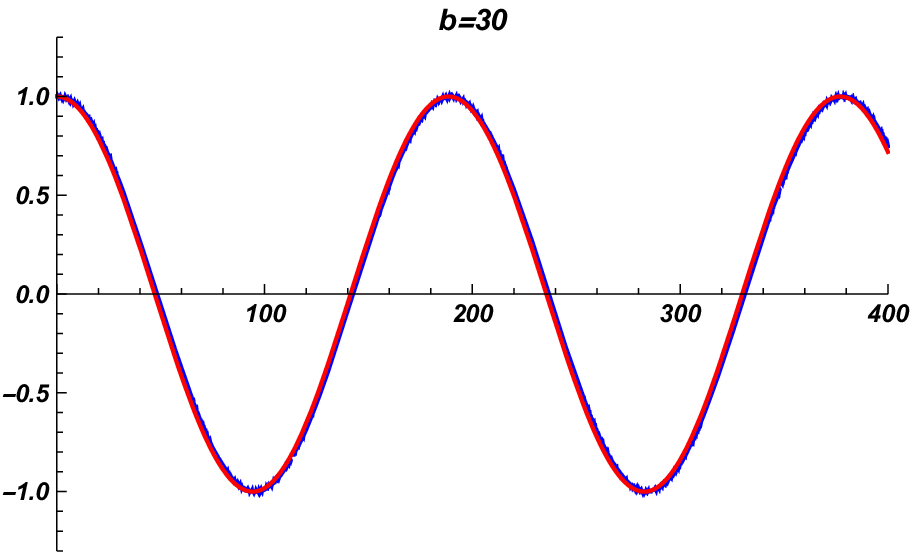}}\\
{\tiny Fig.\ $\!$5$\,\,$ Numerical $u$ (blue) vs asymptotic $u$ (red)
for $\epsilon=1$
with different values of $b$ (and different time-scales).}
\end{center}
\end{figure}}

\def\FigF{\begin{figure}[htb]\begin{center}
\includegraphics[width=4cm]{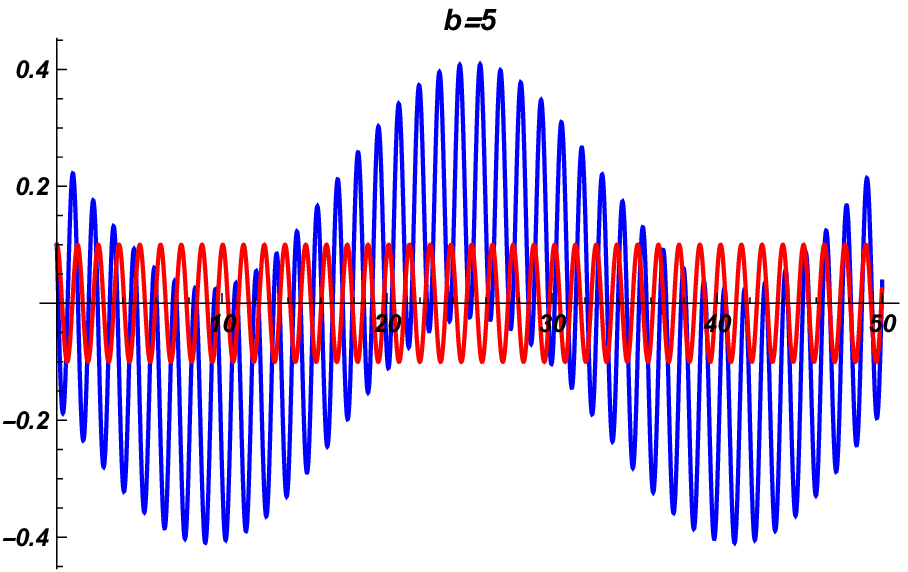}\,\,\,\includegraphics[width=4cm]{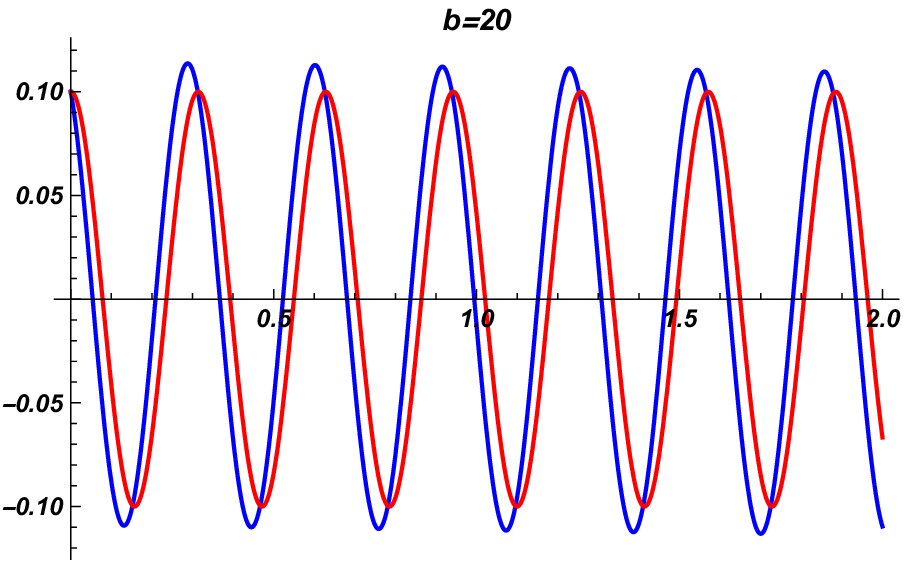}\,\,\,{\includegraphics[width=4cm]{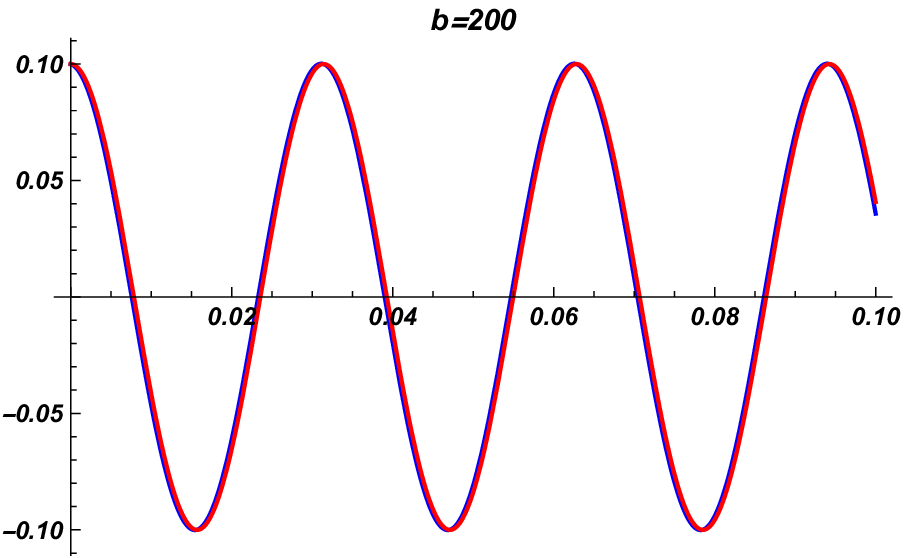}}\\
{\tiny Fig.\ $\!$6$\,\,$ Numerical $v$ (blue) vs asymptotic $v$ (red)
for $\epsilon=1$
with different values of $b$ (and different time-scales).}
\end{center}
\end{figure}}

\def\FigG{\begin{figure}[htb]\begin{center}
\includegraphics[width=7cm]{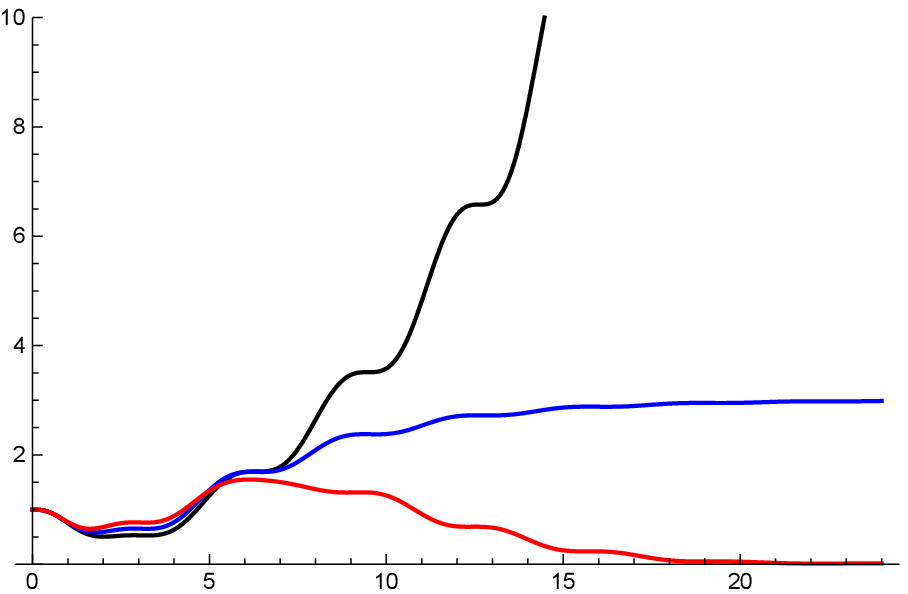}\\
{\tiny Fig.\ $\!$7$\,\,$ Plot of $\E$ for $\epsilon=0.5$ and $b=\sqrt{0.5}-0.1$ (black),
$b=\sqrt{0.5}$ (blue) and $b=\sqrt{0.5}+0.1$ (red).}
\end{center}
\end{figure}}

\def\FigH{\begin{figure}[htb]\begin{center}
\includegraphics[width=7cm]{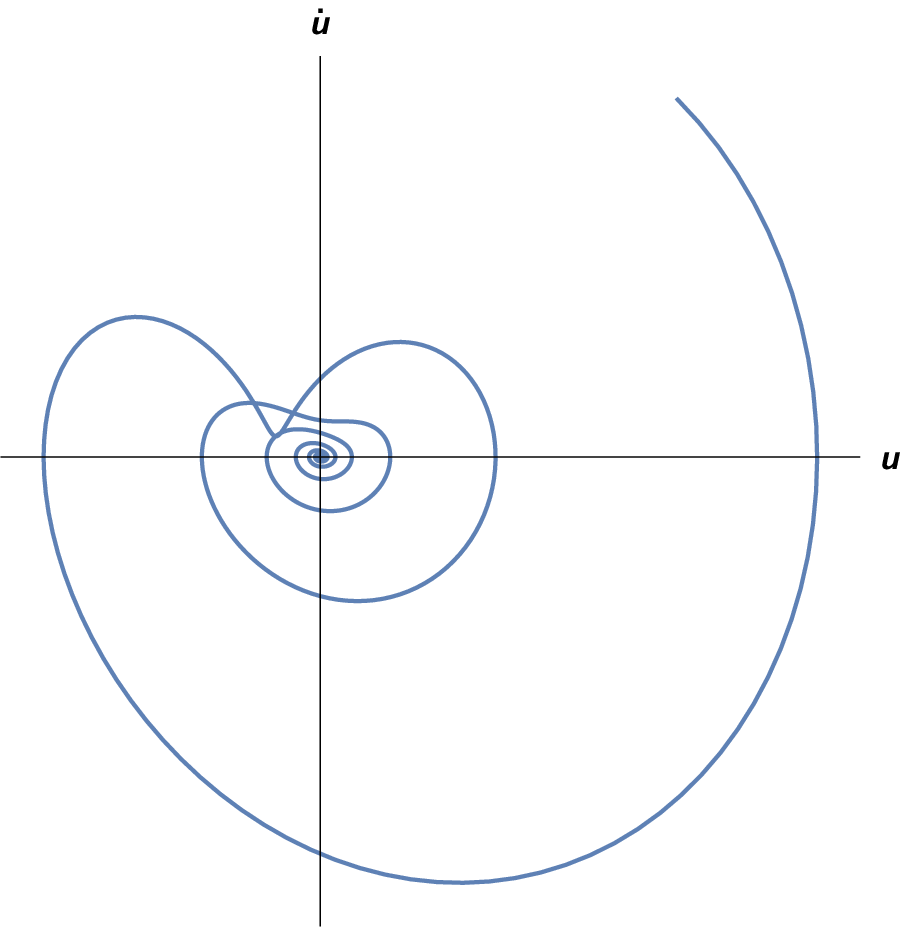}\\
{\tiny Fig.\ $\!$8$\,\,$ Parametric plot of $t\mapsto (u(t),\dot  u(t))$
for $\epsilon=\frac12$ and $b=1$.}
\end{center}
\end{figure}}

\def\FigI{\begin{figure}[htb]\begin{center}
\includegraphics[width=7cm]{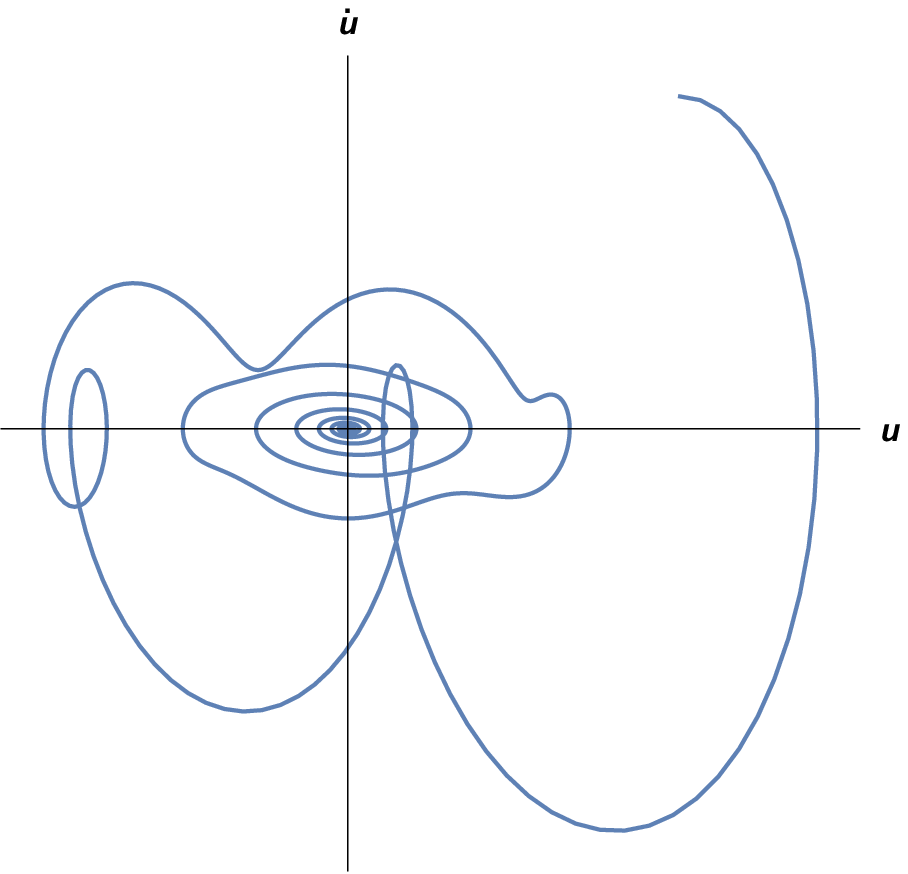}\\
{\tiny Fig.\ $\!$9$\,\,$ Parametric plot of $t\mapsto (u(t),\dot  u(t))$
for $\epsilon=\frac12$ and $b=2$.}
\end{center}
\end{figure}}

\def \no#1#2#3 {{\bf #1} (#3), #2.}
\def \eds#1#2#3 {#1, #2, #3.}

\title[Energy transfer in coupled systems]{A note on the energy transfer\\ in coupled differential systems}

\author[Monica Conti, Lorenzo Liverani and Vittorino Pata]{}

\subjclass{Primary: 34A30, 34D05; Secondary: 35B40, 35L05.}
\keywords{Damped and antidamped equations, coupling parameter, energy transfer, exponential blow up, exponential decay.}

\email{monica.conti@polimi.it}
\email{lorenzo.liverani@polimi.it}
\email{vittorino.pata@polimi.it}

\thanks{$^*$ Corresponding author}

\begin{document}

\maketitle

\centerline{\scshape Monica Conti, Lorenzo Liverani and Vittorino Pata$^*$}
\medskip
{\footnotesize
 \centerline{Politecnico di Milano - Dipartimento di Matematica}
   \centerline{Via Bonardi 9, 20133 Milano, Italy}}

\bigskip

\centerline{(Communicated by the associate editor name)}

\begin{abstract}
We study the energy transfer in the linear system
$$
\begin{cases}
\ddot u+u+\dot u=b\dot v\\
\ddot v+v-\epsilon \dot v=-b\dot u
\end{cases}
$$
made by two coupled differential equations, the first one
dissipative and the second one antidissipative. We see how the competition
between the damping and the antidamping mechanisms
affect the whole system, depending on the coupling parameter $b$.
\end{abstract}

\section{Introduction}

\noindent
The purpose of this work is
to better understand the mutual interaction of two coupled equations,
in terms of the behavior of the associated energy.
What one typically finds in the literature is a system of (ordinary or partial) differential equations, one
of which is conservative and the other one dissipative. The coupling allows
the transfer of dissipation, so that the system becomes globally stable
as time tends to infinity.
Just to quote some results in this direction, we mention the papers \cite{1,2,3,4,5,6,7,8,9,10,11,12,13,14}
and the book~\cite{LZ}, but the list
is far from being exhaustive.

Perhaps, the simplest example is given by
an ideal oscillator without damping,
coupled by velocities
with a physical oscillator subject to dynamical friction,
with initial conditions assigned at time $t=0$.
This is a system of two second-order ODEs of the form
\begin{equation}
\label{ZZZ}
\begin{cases}
\ddot u+u+\dot u=b\dot v,\\
\ddot v+v=-b\dot u,
\end{cases}
\end{equation}
where $b>0$ is the coupling constant.
One is interested to study the longtime behavior
of the associated energy $\E=\E(t)$ given by
$$\E=\frac12\big[u^2+\dot u^2+v^2+\dot v^2\big].$$
Although if $b=0$ the energy of the second equation
is conserved, the effect of the coupling is able to drive $\E(t)$
to zero exponentially fast as $t\to\infty$, no matter how small is $b$.
This result is well known, and can also be obtained as a byproduct of the forthcoming analysis.

Here, instead, we focus on a quite different issue: namely, we want to analyze the effect of the coupling between
a \emph{dissipative} oscillator and an \emph{antidissipative} one.
To this end, we address a simple (yet not so simple) model: namely,
we consider for
$\epsilon>0$ and $b> 0$ the system
\begin{equation}
\label{AAA}
\begin{cases}
\ddot u+u+\dot u=b\dot v,\\
\ddot v+v-\epsilon \dot v=-b\dot u.
\end{cases}
\end{equation}
The situation now is much more intriguing, as we have a competition between
an equation whose solutions decay exponentially fast (in absence of the coupling),
and an equation whose solutions (except the trivial one) exhibit an exponential blow up.

Introducing the four-component (column) vector $\boldsymbol{z}=(u,x,v,y)$, system \eqref{AAA} turns into
the ODE in $\R^4$
\begin{equation}
\label{BBB}
\boldsymbol{\dot z}=\A \boldsymbol{z},
\end{equation}
where the $(4\times 4)$-matrix $\A$ reads
$${\mathbb A}=
\begin{pmatrix*}[r]
0&1&0&0\\
-1&-1&0&b\\
0&0&0&1\\
0&-b&-1&\epsilon
\end{pmatrix*}.
$$
System \eqref{BBB} generates a uniformly continuous semigroup
$$S(t)=e^{t\A},$$
acting by the rule
$$S(t)\boldsymbol{z}_0=\boldsymbol{z}(t),$$
where $\boldsymbol{z}(t)$ is the solution to \eqref{BBB} at time $t$, subject to the initial condition
$\boldsymbol{z}(0)=\boldsymbol{z}_0$. In particular, the energy
corresponding to the initial datum $\boldsymbol{z}_0\in\R^4$ reads
$$\E(t)=\frac12\|S(t) \boldsymbol{z}_0\|^2.$$
Moreover, the
asymptotic properties of $S(t)$ are fully described by the eigenvalues $\lambda_i$
of the matrix $\A$. Indeed,
recalling that the \emph{growth bound} $\omega_\star$ of the semigroup
is defined as
$$\omega_\star=\inf\big\{\omega\in\R\,:\,\|S(t)\|\leq Me^{\omega t}\big\},
$$
for some $M=M(\omega)\geq 1$,
we have the equality
$$\omega_\star=\max_{i} \Re\lambda_i.$$
Here, $\|S(t)\|$ denotes the operator norm of $S(t)$, that is,
$$\|S(t)\|=\sup_{\|\boldsymbol{z}_0\|=1}\|S(t)\boldsymbol{z}_0\|.$$
In particular, when $\omega_\star=0$ the semigroup is bounded (i.e., the energy is bounded for any
initial datum) if and only if all the eigenvalues with null real part are regular. Otherwise,
the norm of $S(t)$ exhibits a blow up of polynomial rate $d$,
where $d$ is the maximum of the defects of those eigenvalues. We address the reader to any classical
ODE textbook for more details (e.g., \cite{HS,PER}).

In summary, the problem reduces to finding such $\lambda_i$, which are the roots of the fourth-order equation
\begin{equation}
\label{poly}
\lambda^4 +(1-\epsilon)\lambda^3  + (2 + b^2 - \epsilon) \lambda^2
+(1-\epsilon)\lambda +1=0.
\end{equation}
Unfortunately, equation \eqref{poly} is not so simple to handle, and the analysis requires some work.

\begin{remark}
Note that the characteristic equation \eqref{poly} depends only on $b^2$. Hence, although we assumed for simplicity
$b>0$, all the subsequent results hold with $b\neq 0$, just replacing every occurrence of $b$ with $|b|$.
\end{remark}

\section{Description of the results}

\noindent
Before entering into technical details, let us anticipate what happens.
When $b$ is small, the two equations do not quite communicate. The result is that the explosive character
of the second equation is predominant, pushing the energy to infinity exponentially fast
for certain initial data.
The two equations start to share the respective energies when $b$ overcomes a certain critical threshold, precisely,
when
$$b>\sqrt{\epsilon}.$$
At this point, the picture strongly depends on the antidamping parameter $\epsilon$.

\begin{itemize}
\item[$\diamond$]
When $\epsilon<1$, the dissipation is stronger than the antidissipation, and the global energy $\E(t)$
undergoes an
exponential decay. The best decay rate is obtained in correspondence of
$$b=\frac{1+\epsilon}2.$$

\item[$\diamond$] On the contrary, when $\epsilon>1$ the dissipation is not enough to contrast
the antidissipation, and the result is an energy which is (generally)
exponentially blowing up for all possible values of $b$.

\smallskip
\item[$\diamond$] The limiting situation is when $\epsilon=1$,
as in that case the damping and the antidamping perfectly compensate. Here,
the system is not conservative,
but nonetheless the energy remains bounded.
Besides, when $b\to \infty$, the energy $\E(t)$ turns into the sum of a highly oscillating term,
possibly vanishing for some particular initial values, and a sinusoid
with a period tending to infinity as well.
\end{itemize}

\section{A detailed discussion}

\noindent
We now proceed to analyze more deeply the three cases.

\begin{aww}
In what follows, for any $z\in{\mathbb C}$, the symbol $\sqrt{z}$
will always mean the value of the complex square root of $z$ whose argument belongs to $(-\frac\pi2,\frac\pi2]$.
With this choice, for any $\alpha,\beta\in\R$ we have
$$\Re\,\sqrt{2(\alpha \pm i \beta)\,}= \sqrt{\rho + \alpha}
\qquad\text{where}\qquad \rho=\sqrt{\alpha^2+\beta^2}.
$$
\end{aww}

Calling now
$$
a= \sqrt{(1+\epsilon)^2-4 b^2},
$$
the four complex roots $\lambda_i$ of equation \eqref{poly} read:
\begin{align*}
\lambda_1 &=\frac{1}{4} \left( \epsilon - 1 + a + \sqrt{2\left( -7 + \epsilon^2 - 2b^2 - (1- \epsilon) a\right)} \right), \\
\lambda_2 &=\frac{1}{4} \left( \epsilon - 1 - a + \sqrt{2\left( -7 + \epsilon^2 - 2b^2 + (1-\epsilon) a\right)} \right), \\
\lambda_3 &=\frac{1}{4} \left( \epsilon - 1 + a - \sqrt{2\left( -7 + \epsilon^2 - 2b^2 - (1- \epsilon) a\right)} \right),\\
\lambda_4 &= \frac{1}{4} \left( \epsilon - 1 - a - \sqrt{2\left( -7 + \epsilon^2 - 2b^2 + (1-\epsilon) a\right)} \right).
\end{align*}

\subsection*{I. The case $\boldsymbol{\epsilon>1}$}
For every value of the coupling parameter $b$, we have
$\omega_\star>0$, meaning that the norm $\|S(t)\|$ of the semigroup
blows up exponentially fast as $t\to\infty$. To show that, it is enough
checking that at least one of the four eigenvalues has
positive real part.
Indeed, since by our convention the square roots have always nonnegative real parts, we readily get
$$
\Re \lambda_1\geq \frac{\epsilon - 1}4 > 0.
$$

\subsection*{II. The case $\boldsymbol{\epsilon=1}$}
Here the four eigenvalues simplify into
\begin{align*}
\lambda_1 &= \frac{1}{2} \left( \sqrt{1- b^2} +\sqrt{-3 - b^2} \right),&\hskip-2.5cm\lambda_4=-\lambda_1,\\
\lambda_2 &= \frac{1}{2} \left( - \sqrt{1- b^2} + \sqrt{-3 - b^2} \right),&\hskip-2.5cm\lambda_3=-\lambda_2.
\end{align*}
We shall distinguish three situations:

\medskip
\noindent
$\bullet$ If $b < 1$, then $\sqrt{1- b^2}>0$.
So $\Re\lambda_{1}>0$, telling that
$\omega_\star>0$.

\medskip
\noindent
$\bullet$ If $b=1$, then
$$
\lambda_1 = \lambda_2 = i\qquad\text{and}\qquad  \lambda_3 = \lambda_4 = -i.
$$
Besides, both the eigenvalues $\pm i$ are not regular, hence with defect 1. Accordingly, $\|S(t)\|$
blows up at infinity with polynomial rate $t$. In fact, in this case we can easily write
the explicit solution corresponding to the generic initial datum
$\boldsymbol{z}_0=(u_0,x_0,v_0,y_0)\in \R^4$ as
\begin{align*}
u(t)&=\frac12\big[(2u_0-t u_0+t v_0) \cos t
+(u_0 - v_0 + 2 x_0 - t x_0 + t y_0) \sin t\big],\\
v(t)&=\frac12\big[ (-t u_0 + 2v_0 + t v_0) \cos t+ (u_0 - v_0 - t x_0 + 2y_0 + ty_0) \sin t\big].
\end{align*}

\smallskip
\noindent
$\bullet$ If $b>1$, then we have four distinct, hence regular, purely imaginary eigenvalues.
This means that there is no uniform decay of the energy, although the energy remains bounded.

\medskip
We conclude the analysis of the case $\epsilon=1$
by examining the qualitative behavior of the solutions for large values of $b$.
When $b\to \infty$, we readily get
$$\lambda_1\sim ib,\qquad\lambda_2\sim \frac{i}b,\qquad\lambda_3\sim-\frac{i}b,\qquad\lambda_4\sim-ib.$$
With the aid of Mathematica{\tiny\texttrademark}, one can compute the asymptotic
form of the matrix ${\mathbb U}$ of the eigenvectors, along with its inverse.
Calling ${\mathbb D}$ the diagonal matrix of the eigenvalues, one can determine
explicitly $S(t)$ via the formula
$$S(t)={\mathbb U}\,e^{t{\mathbb D}}\,{\mathbb U}^{-1}.$$
For $b\to \infty$, this yields
$$S(t)\sim
\begin{pmatrix*}
\cos\frac{t}b & 0 & -\sin\frac{t}b & 0\\
0 & \cos bt & 0 & \sin bt\\
\sin\frac{t}b & 0 & \cos\frac{t}b & 0\\
0 & -\sin bt & 0 & \cos bt
\end{pmatrix*}.
$$
Hence, splitting any initial datum $\boldsymbol{z}_0=(u_0,x_0,v_0,y_0)$ into the sum
$$\boldsymbol{z}_0=\boldsymbol{u}_0+\boldsymbol{x}_0,$$
where $\boldsymbol{u}_0=(u_0,0,v_0,0)$ and $\boldsymbol{x}_0=(0,x_0,0,y_0)$, we obtain the solution
$$\boldsymbol{z}(t)=S(t)\boldsymbol{z}_0\sim\boldsymbol{u}(t)+\boldsymbol{x}(t),
$$
having set
$$\boldsymbol{u}(t)=\textstyle\big(u_0 \cos \frac{t}b -v_0\sin\frac{t}b,\,0\,,u_0 \sin\frac{t}b +v_0\cos\frac{t}b,\,0\big)
$$
and
$$\boldsymbol{x}(t)=\textstyle\big(0,\,x_0 \cos bt +y_0\sin bt,\,0,\,-x_0 \sin bt +y_0\cos bt\big).
$$
So we have the sum of the highly oscillating function $\boldsymbol{x}(t)$ and the sinusoidal
function
$\boldsymbol{u}(t)$ of period $2\pi b\to\infty$. Choosing an initial datum with null velocities, namely, taking
$\boldsymbol{x}_0=\boldsymbol{0}$, in the limiting situation $b=\infty$ we boil down to the constant solution
$\boldsymbol{z}(t)=\boldsymbol{u}_0$.

\subsection*{III. The case $\boldsymbol{\epsilon<1}$}
We show that the exponential decay of the energy occurs when $b>\sqrt{\epsilon}$.
We shall distinguish two situations, depending on the value
$$\eta=\frac{1+\epsilon}2>\sqrt{\epsilon}.
$$

\medskip
\noindent
$\bullet$ If $b\leq \eta$, then
$0\leq a \leq 1 + \epsilon$.
In turn,
$$-7 + \epsilon^2 - 2b^2 \pm (1-\epsilon) a\leq -7 + \epsilon^2 +(1-\epsilon)(1+\epsilon)=-6<0,$$
and consequently
$$
\Re \lambda_1=\Re \lambda_3=\frac14({\epsilon - 1}+a)\qquad\text{and}\qquad
\Re \lambda_2=\Re \lambda_4=\frac14({\epsilon - 1}-a).
$$
At this point, it is convenient to further split the analysis into three subcases.

\smallskip
\noindent
\textbf{-} If $b < \sqrt{\epsilon}$, then
$a>1-\epsilon$. Thus
$\Re \lambda_{1}=\Re \lambda_{3}>0$,
implying that $\omega_\star>0$.

\smallskip
\noindent
\textbf{-} If $b = \sqrt{\epsilon}$, then $a=1-\epsilon$. Therefore,
$$
\Re \lambda_{1}=\Re \lambda_{3} = 0
\qquad\text{and}\qquad \Re\lambda_{2}=\Re\lambda_{4} = \frac{1}{2} \left( \epsilon - 1 \right)<0.
$$
Besides, the four eigenvalues are all distinct, hence regular.
This tells that the energy is bounded, and there exist trajectories
not decaying to zero.

\smallskip
\noindent
\textbf{-} If $\sqrt{\epsilon} < b \leq \eta$, then
$a < 1 - \epsilon$, which immediately gives
$\Re \lambda_i < 0$ for all $i$. The energy undergoes an exponential decay.

\medskip
\noindent
$\bullet$ If $b>\eta$, then
$$a=i\sqrt{4b^2-(\epsilon+1)^2}.$$
Therefore,
$$
\Re\,\sqrt{2\left( -7 + \epsilon^2 - 2b^2 \pm (1-\epsilon)a\right)}= \sqrt{\rho -7 + \epsilon^2 - 2b^2},
$$
where
$$
\rho = 2\sqrt{b^4 +2b^2(4-\epsilon)+3(4-\epsilon^2)}.
$$
Accordingly,
\begin{align*}
\label{span}
\Re \lambda_1 & =\Re\lambda_2 = \frac14\left(\epsilon - 1 + \sqrt{\rho -7 + \epsilon^2 - 2b^2}\right),\\
\Re \lambda_3 & =\Re\lambda_4 = \frac14\left(\epsilon - 1 - \sqrt{\rho -7 + \epsilon^2 - 2b^2}\right).
\end{align*}
It is then clear that $\Re\lambda_3=\Re\lambda_4<0$, and with standard computations we readily check that
$\Re\lambda_1 = \Re \lambda_2 < 0$ as well.

\subsection*{IV. The best decay rate}
Once we know that when $\epsilon<1$ and $b>\sqrt{\epsilon}$ the exponential decay occurs,
it is interesting to establish for which value of the coupling parameter $b$ the best decay rate
is attained. From the previous discussion, it is readily seen that when $\eta\neq b>\sqrt{\epsilon}$ then
$$\Re\lambda_1>\frac{\epsilon-1}4\qquad\Rightarrow\qquad \omega_\star>\frac{\epsilon-1}4.
$$
Accordingly, the smallest possible value
is exactly
$$\omega_\star=\frac{\epsilon-1}4,$$
which is achieved when $b=\eta$. In this case,
\begin{align*}
\lambda_1=\lambda_2 &= \frac{1}{4} \left( \epsilon - 1 + \sqrt{-15 + \epsilon^2 - 2\epsilon} \right), \\
\lambda_3 =\lambda_4 &= \frac{1}{4} \left( \epsilon - 1 - \sqrt{-15 + \epsilon^2 - 2\epsilon} \right),
\end{align*}
and the two distinct eigenvalues, sharing the same real part, can be shown to be
nonregular, hence with defect 1.
Then the optimal exponential decay rate $(1-\epsilon)/4=-\omega_\star$ for the semigroup norm
is polynomially penalized, yielding the best possible decay estimate
$$\|S(t)\|\leq C(1+t)e^{-\frac{1-\epsilon}4t},$$
for some $C\geq 1$.
Observe also that $\Re \lambda_1=\Re\lambda_2\to 0$ when $b\to\infty$,
telling that the exponential decay rate tends to zero when $b$ becomes large.
Indeed, we have the asymptotic expansion
$\rho=2b^2+8-2\epsilon+o(1)$ as $b\to\infty$.

\medskip
\begin{remark}
The reader will have no difficulty to ascertain
that the analysis made in the previous points III and IV covers
the limit value $\epsilon=0$ as well, corresponding to problem~\eqref{ZZZ}.
Here, $b=\eta=\frac12$, and the two (nonregular) distinct eigenvalues read
$$\frac{1}{4} \left(- 1 \pm i\sqrt{15} \right).$$
The optimal decay estimate becomes
$$\|S(t)\|\leq C(1+t)e^{-\frac14t}.$$
\end{remark}

\section{The infinite dimensional case}

\noindent
The finite-dimensional analysis carried out so far, besides having an interest by itself,
can also be extended to cover some infinite-dimensional models.
Indeed, in greater generality, one might consider the same problem for the
system
\begin{equation}
\label{QQQ}
\begin{cases}
\ddot u+Au+\dot u=b\dot v,\\
\ddot v+Av-\epsilon \dot v=-b\dot u,
\end{cases}
\end{equation}
where $A$ is a strictly positive selfadjoint operator acting on a Hilbert space $H$, with
compactly embedded domain ${\mathfrak D}(A)\Subset H$. From the classical
theory of semigroups~\cite{PAZ}, system \eqref{QQQ} is well known to
generate a strongly continuous semigroup $S(t)$ acting on the product Hilbert space
$${\mathcal H}={\mathfrak D}(A^{\frac12})\times H \times {\mathfrak D}(A^{\frac12})\times H.
$$
A concrete realization of \eqref{QQQ} is the system of PDEs
$$
\begin{cases}
u_{tt}-\Delta u+ u_{t}=b v_{t},\\
v_{tt}-\Delta v-\epsilon v_{t}=-b u_{t},
\end{cases}
$$
where $\Delta$ is the Laplace-Dirichlet operator
acting on the Hilbert space $L^2(\Omega)$, for some  bounded domain $\Omega\subset\R^N$
with boundary $\partial\Omega$ smooth enough.

Here the picture is exactly the same as in the ODE system considered before.
The desired results can be proved by projecting the equations on the eigenvectors
of $A$, and then by computing the decay rate of each single mode.
The only difference occurs in the case $\epsilon<1$, where $b=(1+\epsilon)/2$
is still the value corresponding to the best exponential decay rate, but the decay rate
itself can be affected by the
first eigenvalue $\lambda_1>0$ of $A$.
This happens when $\lambda_1$ is small. Indeed,
if
$\lambda_1\geq (1-\epsilon)^2/{16}$,
then we recover
the exponential decay rate $(1-\epsilon)/4$, up to a polynomial correction.

\begin{remark}
In fact, the request that $A$ has compact inverse is not really needed, although
this assumption greatly simplifies the analysis, since in this case the spectrum of $A$
is made by eigenvalues only. If $A^{-1}$ is not compact, a deeper use of the spectral theory
and the related functional calculus is required. We refer the interested reader to the paper~\cite{GGPM},
where these techniques have been successfully exploited in the analysis of the best
exponential decay rate for an abstract
weakly damped wave equation.
\end{remark}

\section{Some figures}

\noindent
We conclude the paper with some figures illustrating our analysis.
We will concentrate on the two more interesting cases $\epsilon=1$ and $\epsilon<1$.

\medskip
\bigskip
The first set of figures concerns with the case $\epsilon=1$.

\medskip
\noindent $\diamond$ In Fig.\ 1 we see the behavior of the energy $\E(t)$ corresponding to
the initial value $\boldsymbol{z}_0=(1,0,0,0)$, for $b<1$ (exponential blow up),
$b=1$ (polynomial blow up of rate $t^2$), and $b>1$ (bounded energy).

\medskip
\bigskip
\FigA
\bigskip
\bigskip
\bigskip

\medskip
\noindent $\diamond$
In Fig.\ 2, again for the
initial value $\boldsymbol{z}_0=(1,0,0,0)$, we represent the phase portrait of
the first component $u(t)$ of the solution along with its derivative $\dot u (t)$.
If one takes (as in the figure)
$b=\sqrt{q+q^{-1}-1},$
with $q$ rational number, then the phase portrait becomes periodical.

\newpage

\bigskip
\FigB
\bigskip

\medskip
\noindent $\diamond$ In Fig.\ 3 and Fig.\ 4 we plot the energy $\E(t)$
for three values of $b>1$. In Fig.\ 3 the initial value is $\boldsymbol{z}_0=(1,0,0,0)$.
Here, we see that as $b$ increases the energy becomes sinusoidal. Instead,
in Fig.\ 4 we take the initial value $\boldsymbol{z}_0=(1,0.5,0,0)$. We
observe that the oscillations about the sinusoid persist, and their frequency increases dramatically
as $b\to \infty$.

\bigskip
\FigC
\bigskip
\FigD
\bigskip

\newpage
\noindent $\diamond$ In Fig.\ 5 and Fig.\ 6 we compare for different values
of $b$ the numerical solutions $u(t)$ and $v(t)$ with their asymptotic counterparts found in Section~3 part II.
In both cases, we take the initial value $\boldsymbol{z}_0=(1,0.1,0,0)$.
As predicted, the two curves overlap when $b\to \infty$.

\bigskip
\FigE
\bigskip\bigskip
\FigF
\bigskip

\bigskip

Finally, we focus on the case $\epsilon<1$.

\medskip
\noindent $\diamond$ In Fig.\ 7 we plot the energy corresponding to the initial value $\boldsymbol{z}_0=(1,0,0,0)$
for $b<\sqrt{\epsilon}$ (exponential blow up), $b=\sqrt{\epsilon}$ (bounded energy), and $b>\sqrt{\epsilon}$
(exponential decay).

\bigskip
\FigG
\bigskip

\medskip
\noindent $\diamond$
In Fig.\ 8 and Fig.\ 9, taking $\epsilon=\frac12$ and
the initial datum $\boldsymbol{z}_0=(1,1,1,1)$, we represent the phase portrait of
the first component $u(t)$ of the solution along with its derivative $\dot u (t)$ for
$b=1$ and $b=2$, respectively.

\bigskip
\FigH

\bigskip
\FigI

\section*{Acknowledgments} We thank Professor Giulio Magli for fruitful discussion and comments.


\medskip
Received xxxx 20xx; revised xxxx 20xx.
\medskip

\end{document}